\documentclass[11pt,a4paper]{article}
\usepackage[margin=28mm]{geometry}
\usepackage[T1]{fontenc}
\usepackage{lmodern,microtype}
\usepackage{amsmath,amssymb,amsthm,mathtools}
\usepackage{enumitem}
\usepackage[colorlinks=true,linkcolor=blue!50!black,citecolor=blue!50!black,urlcolor=blue!50!black]{hyperref}
\usepackage{xcolor}
\usepackage{authblk}
\newtheorem{theorem}{Theorem}[section]
\newtheorem{lemma}[theorem]{Lemma}

\newtheorem{remark}[theorem]{Remark}

\newtheorem{definition}[theorem]{Definition}
\newcommand{\Om}{\Omega}\newcommand{\Ga}{\Gamma}\newcommand{\G}{\Gamma_1}
\newcommand{\R}{\mathbb R}
\newcommand{\du}[2]{\langle #1,#2\rangle}
\newcommand{\norm}[1]{\lVert #1\rVert}
\newcommand{\eps}{\varepsilon}
\title{Nonsmooth Multivalued Leak Problem for the Evolutionary Two-Dimensional Navier--Stokes Equations}
\author{Paweł Szafraniec}
\affil{Faculty of Mathematics and Computer Science \\
	Jagiellonian University in Krak\'ow \\
	ul. \L{}ojasiewicza 6, 30-348 Krak\'ow, Poland}
\begin{document}
\maketitle
\begin{abstract}
We study the evolutionary incompressible Navier--Stokes equations in two space dimensions subject to a homogeneous no-slip condition on one part of the boundary and a nonsmooth, possibly nonmonotone, leak condition on the remaining part. The leak law is expressed through the Clarke subdifferential of a boundary potential and involves the physical Cauchy normal stress, rather than a Bernoulli pressure. We establish existence of a global weak solution for a class of potentials containing a coercive cubic component and a locally Lipschitz perturbation of at most linear subgradient growth. 
\end{abstract}

\section{Introduction}
Let an incompressible viscous fluid occupy a bounded planar domain whose boundary is divided into a no-slip part and a part through which fluid may enter or leave. In the classical Dirichlet problem the convective term becomes zero, when tested with $u$, i.e. $b(u,u,u)=0$, same for slip problems. This property fails for a leak boundary: for a velocity field satisfying the zero tangential velocity condition, integration by parts gives, unfortunately
\begin{equation}\label{introflux}
 \int_\Om ((u\cdot\nabla)u)\cdot u\,dx
 =\frac12\int_{\G}u_\nu |u|^2\,d\Ga
 =\frac12\int_{\G}u_\nu^3\,d\Ga.
\end{equation}
The sign of the boundary integral is not prescribed. This is the principal obstacle in deriving energy estimates for the leak problem with the physical Cauchy stress.

Variational formulations of slip and leak boundary conditions have been investigated in, among others, \cite{Fujita94,Fujita01,Kashiwabara,AnLiLi}. Nonsmooth and nonconvex boundary laws are naturally formulated using Clarke's generalized gradient, leading to hemivariational inequalities; see \cite{Clarke,MigorskiOchalSofonea,Naniewicz,MigorskiOchal07}. Related evolutionary problems occur in \cite{Szafraniec17}. The issue addressed here is the interaction of the boundary flux in \eqref{introflux} with a nonmonotone leak law.

Usually in this type of problems, a linear-growth hypothesis on the entire Clarke subdifferential allows control of the boundary work by lower-order terms. It does not, by itself, absorb the cubic flux in \eqref{introflux}. A strategy based on assuming the normal trace small for a short time consequently requires a uniform, approximation-independent time interval, which is not supplied by the usual $L^\infty(H)\cap L^2(V)$ estimate, as author realized after many failed attempts. We instead impose a dissipative cubic part in the potential,
\begin{equation}\label{introdecomp}
 j(x,t,r)=\frac\kappa3|r|^3+\psi(x,t,r),\qquad \kappa>\frac12,
\end{equation}
while retaining a genuinely nonsmooth and nonmonotone perturbation $\psi$. This is a substantive constitutive restriction, not a consequence of the original linear-growth assumption. It yields global weak existence on each finite time interval. 

We develop the functional setting in Section~\ref{sec:prelim}, formulate the physical and weak problems in Section~\ref{sec:problem}, prove existence in Section~\ref{sec:proof}, and discuss examples and related formulations in Section~\ref{sec:discussion}.

\section{Functional setting and preliminary results}\label{sec:prelim}
Let $\Om\subset\R^2$ be a bounded domain with $C^2$ boundary $\Gamma$ consiting of two parts $\Gamma_0$ and $\Gamma_1$ whose closures cover $\Gamma$, i.e. $\Gamma = \bar{\Gamma_0} \cup \bar{\Gamma_1}$ with $meas(\Gamma_0)>0$. The outer unit normal vector is denoted by $\nu$. For a vector field $v$ on the boundary we decompose it into normal and tangential components, as follows:
\[
 v_\nu=v\cdot\nu,\qquad v_\tau=v-v_\nu\nu.
\]
We will work with the natural real Hilbert spaces
\begin{align*}
 V&=\{v\in H^1(\Om;\R^2):\operatorname{div}v=0,\ v=0\text{ on }\Ga_0,\ v_\tau=0\text{ on }\G\},\\
 H&=\overline{V}^{\,L^2(\Om;\R^2)}.
\end{align*}
 Identifying $H$ with its dual gives the Gelfand triple $V\hookrightarrow H\hookrightarrow V^*$. We use the $L^2$ scalar product $(\cdot,\cdot)_H$.

Write $\eps(v)=\tfrac12(\nabla v+\nabla v^\top)$ for a linearised stress tensor and introduce binlinear form $a_0$ and trilinear form $a_1$ respecively by
\begin{align}
 a_0(u,v)&=2\nu_0\int_\Om\eps(u):\eps(v)\,dx,\qquad \nu_0>0,\label{a0}\\
 a_1(u,v,w)&=\int_\Om ((u\cdot\nabla)v)\cdot w\,dx.\label{a1}
\end{align}
Korn's inequality, using $meas(\Gamma_0)>0$, ensures the coercivity
\begin{equation}\label{coerc}
 a_0(v,v)\ge \alpha\norm{v}_V^2,\qquad v\in V,
\end{equation}
where $\norm{\cdot}_V$ is an equivalent $H^1$ norm and $\alpha>0$.

In two dimensions we have $V\hookrightarrow L^4(\Om;\R^2)$ and
\begin{equation}\label{convbound}
 |a_1(u,v,w)|\le C\norm u_{L^4(\Om)}\norm v_V\norm w_{L^4(\Om)}\le C\norm u_V\norm v_V\norm w_V.
\end{equation}
For $u,v\in V$, the divergence theorem gives, unfortunetly for us:
\begin{equation}\label{boundaryid}
 a_1(u,v,v)=\frac12\int_{\G}u_\nu|v|^2\,d\Ga.
\end{equation}
For smooth fields, the product rule gives
\[((u\cdot\nabla)v)\cdot v
=\tfrac12\operatorname{div}(u|v|^2)
-\tfrac12(\operatorname{div}u)|v|^2.\]
Integrate over $\Om$. The divergence term in the volume vanishes,
while the boundary integral vanishes on $\Ga_0$ because $u=0$.
On $\G$, $v=v_\nu\nu$, and therefore $|v|^2=v_\nu^2$.
The trace embedding $V\hookrightarrow L^4(\G)$ and density extend
the identity to all fields in $V$.

In particular, because $v_\tau=0$ on $\G$, setting $u=v$ gives \eqref{introflux}.

\

We also have the trace theorem: For each $\delta>0$ there exists $C_\delta$ such that
\[\|v_\nu\|_{L^2(\G)}^2\le\delta\|v\|_V^2+C_\delta\|v\|_H^2\quad(v\in V).\]

We recall the Clarke generalized directional derivative of a locally Lipschitz function $\varphi:\R\to\R$,
\[
 \varphi^0(r;s)=\limsup_{y\to r,\,h\downarrow0}\frac{\varphi(y+hs)-\varphi(y)}h,
\]
and the Clarke subdifferential
\[
 \partial\varphi(r)=\{\zeta\in\R:\zeta s\le\varphi^0(r;s)\ \text{for all }s\in\R\}.
\]
 If $\varphi_1$ is continuously differentiable and $\varphi_2$ is locally Lipschitz, then
\begin{equation}\label{sumrule}
 \partial(\varphi_1+\varphi_2)(r)=\varphi_1'(r)+\partial\varphi_2(r).
\end{equation}
We shall use the compactness criterion for functions bounded in $L^2(0,T;V)$ whose time derivatives are bounded in $L^1(0,T;V^*)$; see \cite{Simon}. We emphasize that the exponent $1$ for the time derivative is sufficient for compactness in $L^2(0,T;H)$, although it does not imply the usual Hilbert-space evolution regularity $W(0,T)$.
We also give two calssical lemmas used throughout

\begin{lemma}[Trace theorem]\label{lem:trace}
	Suppose $(w_m)$ is bounded in $L^\infty(0,T;H)\cap L^2(0,T;V)$ and $(w_m')$ is bounded in $L^1(0,T;V^*)$. Then $(\gamma_\nu w_m)$ is relatively compact in $L^2(0,T;L^2(\G))$.
\end{lemma}

\begin{proof}
	Let $s\in(1/2,1)$. Taking into account the embeddings
	\[
	V\subset H^s(\Omega;\mathbb{R}^2)
	\hookrightarrow V^*,
	\]
	where the first embedding is compact,  and using the boundedness of the sequences $(w_m)$
	in $L^2(0,T;V)$ and $(w_m')$ in $L^1(0,T;V^*)$,
	from \cite[Corollary 4, p.~85]{Simon} we deduce
	that $(w_m)$ is relatively compact in
	$L^2(0,T;H^s(\Omega;\mathbb{R}^2))$.
	The assertion now follows from the continuity of
	the normal trace operator.
	\[
	\gamma_\nu:H^s(\Omega;\mathbb{R}^2)
	\to L^2(\Gamma_1).
	\]
\end{proof}

\begin{lemma}[Erhling Lemma]\label{erhling}
	For each $\varepsilon>0$ there exists $C_\varepsilon$ such that
	\[\|v_\nu\|_{L^2(\G)}^2\le\varepsilon\|v\|_V^2+C_\varepsilon\|v\|_H^2, \quad \quad  v\in V.\]
\end{lemma}

\section{The leak problem and hypotheses}\label{sec:problem}
For a finite $T>0$, let $Q=\Om\times(0,T)$ and $\Sigma_i=\Ga_i\times(0,T)$. We consider
\begin{align}
 \partial_tu-\nu_0\Delta u+(u\cdot\nabla)u+\nabla p&=f&&\text{in }Q,\label{NS}\\
 \operatorname{div}u&=0&&\text{in }Q,\label{div}\\
 u&=0&&\text{on }\Sigma_0,\label{noslip}\\
 u_\tau=0,\qquad -\sigma_\nu(u,p)&\in\partial_rj(x,t,u_\nu)&&\text{on }\Sigma_1,\label{leak}\\
 u(0)&=u_0&&\text{in }\Om.\label{initial}
\end{align}
Here $\sigma(u,p)=-pI+2\nu_0\eps(u)$ is the Cauchy stress tensor and $\sigma_\nu=(\sigma\nu)\cdot\nu$. The law \eqref{leak} involves the physical normal stress and not a pressure modified by the kinetic energy, as it is usually considered in the litereture. Main idea of this paper is to show existence in this more natural setting.

The assumptions are as follows.
\begin{description}[leftmargin=14mm,style=nextline]
\item[\textbf{$H_0$.}] $f\in L^2(0,T;V^*)$, $u_0\in H$, and $\nu_0>0$. 
\item[\textbf{$H_j^{\rm cub}$.}] There exist $\kappa>\frac12$ and $M>0$ such that
\begin{equation}\label{hypj}
 j(x,t,r)=\frac\kappa3|r|^3+\psi(x,t,r),\qquad
 |\zeta|\le M(1+|r|)\quad(\zeta\in\partial_r\psi(x,t,r)).
\end{equation}
For almost every $(x,t)\in\Sigma_1$, the function $r\mapsto\psi(x,t,r)$ is locally Lipschitz. For every $r\in\R$, $(x,t)\mapsto\psi(x,t,r)$ is measurable, and $\psi(\cdot,\cdot,0)\in L^1(\Sigma_1)$. We assume the associated Clarke multifunction is measurable in $(x,t)$; this holds, for example, for the usual measurable locally Lipschitz normal-integrand setting.
\end{description}
The growth bound in \eqref{hypj} is imposed on $\partial\psi$, not on $\partial j$. By \eqref{sumrule},
\begin{equation}\label{decomposition}
 \partial_rj(x,t,r)=\kappa|r|r+\partial_r\psi(x,t,r).
\end{equation}
The perturbation $\psi$ is therefore allowed to be nonconvex and nonregular.

\begin{definition}[Weak solution]\label{def:weak}
A weak solution on $(0,T)$ is a pair $(u,\xi)$ satisfying
\begin{align*}
 &u\in L^\infty(0,T;H)\cap L^2(0,T;V),\quad u_\nu\in L^3(\Sigma_1),\\
 &\xi\in L^2(\Sigma_1),\quad \xi(x,t)\in\partial_r\psi(x,t,u_\nu(x,t))\quad\text{a.e.},\\
 &u'\in L^1(0,T;V^*),\quad \eta:=\kappa|u_\nu|u_\nu+\xi\in L^{3/2}(\Sigma_1),
\end{align*}
such that for every $v\in V$ we have the equality
\begin{equation}\label{variational}
 \du{u'(t)}v+a_0(u(t),v)+a_1(u(t),u(t),v)
 +\int_{\G}\eta(x,t)v_\nu(x)\,d\Ga=\du{f(t)}v
\end{equation}
for almost every $t\in(0,T)$. The initial condition is $u(0)=u_0$ in the weak $H$ sense.
\end{definition}
The terms in \eqref{variational} are integrable in time: the convection operator belongs to $L^1(0,T;V^*)$, while the cubic boundary operator belongs to $L^{3/2}(0,T;V^*)$. Thus the natural time derivative is $L^1(V^*)$.

\begin{theorem}[Global weak existence]\label{thm:main}
Assume $H_0$ and $H_j^{\rm cub}$. Then for every finite $T>0$ problem \eqref{variational} admits a weak solution in the sense of Definition~\ref{def:weak}. There is a constant $C_T$, depending on $T$, $\Om$, $\nu_0$, $\kappa$, $M$, $\norm{u_0}_H$, and $\norm f_{L^2(0,T;V^*)}$, such that
\begin{equation}\label{mainenergy}
 \norm u_{L^\infty(0,T;H)}^2+\norm u_{L^2(0,T;V)}^2
 +(\kappa-\tfrac12)\norm{u_\nu}_{L^3(\Sigma_1)}^3\le C_T.
\end{equation}
No uniqueness assertion ia claimed, as that would force stronger hypothesis on $j$.
\end{theorem}

\section{Proof of Theorem~\ref{thm:main}}\label{sec:proof}
\subsection{Spectral Galerkin approximation and regularization}
The form $a_0$ is symmetric, continuous and coercive on $V$, and $V$ is compactly embedded in $H$. Hence the operator connected to it is positive self-adjoint in $H$, so it has an $H$-orthonormal basis made of eigenvectors $(e_k)_{k\ge1}\subset V$, orthogonal also with respect to $a_0$. Let
\[
 V_m=\operatorname{span}\{e_1,\ldots,e_m\},\qquad P_m:H\to V_m
\]
be the $H$-orthogonal projection. It extends by duality to $V^*$ and satisfies
\begin{equation}\label{proj}
 \norm{P_m v}_V\le\norm v_V,\qquad
 \norm{P_mF}_{V^*}\le\norm F_{V^*},\qquad v\in V,\ F\in V^*.
\end{equation}
The second assertion follows by testing $P_mF$ against $v\in V$ and using the first. More explicitly, if $a_0(e_k,v)=\lambda_k(e_k,v)_H$, then
$\norm{P_mv}_V^2=\sum_{k\le m}\lambda_k|(v,e_k)_H|^2
\le\sum_{k\ge1}\lambda_k|(v,e_k)_H|^2=\norm v_V^2$
when $\norm v_V^2=a_0(v,v)$; replacing the norm by an equivalent one
only changes the projection constant. This uniform stability permits the
finite-dimensional equations to be interpreted in the fixed space $V^*$.
 Put $u_{0m}=P_mu_0$; then $u_{0m}\to u_0$ in $H$ and $\norm{u_{0m}}_H\le\norm{u_0}_H$.

Choose a nonnegative smooth mollifier $\rho$ supported in $[-1,1]$ with unit integral and define $\rho_m(s)=m\rho(ms)$. Set
\[
 \psi_m(x,t,r)=\int_\R\rho_m(s)\psi(x,t,r-s)\,ds,\qquad
 g_m(x,t,r)=\partial_r\psi_m(x,t,r).
\]
The scalar functions $\psi_m(x,t,\cdot)$ are smooth. The growth condition gives, uniformly in $m\ge1$,
\begin{equation}\label{gm}
 |g_m(x,t,r)|\le M(1+|r|+m^{-1})\le 2M(1+|r|).
\end{equation}
For almost every $(x,t)$, $g_m(x,t,r)$ belongs to the closed convex hull of the values of $\partial_r\psi(x,t,\cdot)$ in the interval $[r-m^{-1},r+m^{-1}]$.

We seek an absolutely continuous $V_m$-valued function $u_m$ satisfying, for $v\in V_m$ and almost every $t$,
\begin{align}\label{galerkin}
 (u_m',v)_H+a_0(u_m,v)+a_1(u_m,u_m,v)
 &+\kappa\int_{\G}|u_{m\nu}|u_{m\nu}v_\nu\,d\Ga\\
 &+\int_{\G}g_m(x,t,u_{m\nu})v_\nu\,d\Ga=\du{f(t)}v,\nonumber
\end{align}
with $u_m(0)=u_{0m}$. The right-hand side of the corresponding finite-dimensional ODE is measurable in $t$, continuous in the coefficients, and locally bounded by an integrable function. Carath\'eodory's theorem gives a local solution.

Now take $v=u_m(t)$ in \eqref{galerkin}. The identity \eqref{introflux} gives
\begin{align}\label{energy1}
 \frac12\frac{d}{dt}\norm{u_m}_H^2+\alpha\norm{u_m}_V^2
 +\kappa\int_{\G}|u_{m\nu}|^3\,d\Ga
 &=\du f{u_m}-\frac12\int_{\G}u_{m\nu}^3\,d\Ga\\
 &\quad-\int_{\G}g_m(x,t,u_{m\nu})u_{m\nu}\,d\Ga.\nonumber
\end{align}
The convective flux has no fixed sign, but still we can conclude that
\begin{equation}\label{absorb}
 \kappa|r|^3+\frac12r^3\ge(\kappa-\tfrac12)|r|^3\qquad(r\in\R).
\end{equation}
By \eqref{gm} we have
\begin{equation}\label{perturb}
 \left|\int_{\G}g_m(x,t,u_{m\nu})u_{m\nu}\,d\Ga\right|
 \le C\bigl(1+\norm{u_{m\nu}}_{L^2(\G)}^2\bigr).
\end{equation}
Combining Lemma~\eqref{erhling} and \eqref{perturb}, and using Young's inequality gives us
\begin{equation}\label{energy2}
 \frac{d}{dt}\norm{u_m}_H^2+c_0\norm{u_m}_V^2
 +2(\kappa-\tfrac12)\norm{u_{m\nu}}_{L^3(\G)}^3
 \le C\bigl(1+\norm{f(t)}_{V^*}^2+\norm{u_m}_H^2\bigr)
\end{equation}
with $c_0>0$ independent of $m$. For clarity, choose $\varepsilon_1>0$ in the force estimate
$|\langle f,u_m\rangle|\le\varepsilon_1\norm{u_m}_V^2+
(4\varepsilon_1)^{-1}\norm f_{V^*}^2$.
Next, in the trace estimate \eqref{erhling} choose $\varepsilon_2>0$
so that the contribution of the perturbation to $\norm{u_m}_V^2$
is less than $\alpha/4$. Choosing $\varepsilon_1$ in the same way
leaves a positive multiple of the viscous energy. The constants
$\varepsilon_1,\varepsilon_2$ are fixed independently of $m$ and of
the length of the interval. In particular, no bound on
$\norm{u_{m\nu}(t)}_{L^2(\G)}$ is used to obtain coercivity.
Gronwall's lemma and integration gives
\begin{equation}\label{apriori}
 \sup_{0\le t\le T}\norm{u_m(t)}_H^2
 +\int_0^T\norm{u_m}_V^2\,dt
 +\int_0^T\!\int_{\G}|u_{m\nu}|^3\,d\Ga\,dt\le C_T.
\end{equation}
In particular, the finite-dimensional solution cannot blow up on $[0,T]$. 

We now test the Galerkin equation with $v=u_m$, and it becomes
\begin{align*}
\frac12 \frac{d}{dt} \|u_m\|_H^2+a_0(u_m,u_m)
 +\int_{\G}\left(\kappa|u_{n\nu}|^3+\frac12 u_{n\nu}^3\right)d\Ga
 =\langle f,u_m\rangle-\int_{\G} g_m(u_{n\nu})u_{n\nu}\,d\Ga.
\end{align*}
For $r\ge0$ the boundary integrand equals
$(\kappa+\frac12)|r|^3$, whereas for $r<0$ it equals
$(\kappa-\frac12)|r|^3$. Consequently,
\begin{align}\label{extraE}
\frac12 \frac{d}{dt} \|u_m\|_H^2+\alpha\|u_m\|_V^2
 +(\kappa-\tfrac12)\|u_{n\nu}\|_{L^3(\G)}^3
 \le |\langle f,u_m\rangle|
 +\left|\int_{\G} g_m(u_{n\nu})u_{n\nu}\,d\Ga\right|.
\end{align}
Choose $\varepsilon_1>0$. Young's inequality gives
\[|\langle f,u_m\rangle|
 \le\varepsilon_1\|u_m\|_V^2
 +(4\varepsilon_1)^{-1}\|f\|_{V^*}^2.\]
Using $|g_m(r)|\le2M(1+|r|)$ and $2|r|\le1+r^2$,
\begin{align*}
\left|\int_{\G} g_m(u_{m\nu})u_{m\nu}\,d\Ga\right|
 &\le2M\int_{\G}(|u_{m\nu}|+|u_{m\nu}|^2)\,d\Ga\\
 &\le M|\G|+3M\|u_{m\nu}\|_{L^2(\G)}^2\\
 &\le M|\G|+3M\varepsilon_2\|u_m\|_V^2
 +3MC_{\varepsilon_2}E_m.
\end{align*}
Choose $\varepsilon_1=\alpha/4$ and
$\varepsilon_2=\alpha/(12M)$ when $M>0$ (the case $M=0$
is obvious). Then $\varepsilon_1+3M\varepsilon_2=\alpha/2$.
Multiplying \eqref{extraE} by $2$ gives
\begin{equation}\label{extraG}
 \frac{d}{dt} \|u_m(t)\|_H^2+\alpha\|u_m\|_V^2
 +2(\kappa-\tfrac12)\|u_{m\nu}\|_{L^3(\G)}^3
 \le C_0+C_1E_m+C_2\|f\|_{V^*}^2.
\end{equation}
Multiplying \eqref{extraG} by $e^{-C_1t}$ and omiting its two
nonnegative terms and integrating from $0$ to $t$ gives
\begin{align*}
 \|u_m(t)\|_H^2&\le e^{C_1t}\|u_m(0)\|_H^2
 +\int_0^t e^{C_1(t-s)}(C_0+C_2\|f(s)\|_{V^*}^2)\,ds\\
 &\le e^{C_1T}\bigl(\|u_0\|_H^2+C_0T
 +C_2\|f\|_{L^2(0,T;V^*)}^2\bigr)=:K_T.
\end{align*}
Integrating \eqref{extraG} and using the fact that $\|u_m(T)\|_H\ge0$ we find that
\begin{align*}
\alpha\int_0^T\|u_m\|_V^2dt
 +2(\kappa-\tfrac12)\int_0^T\|u_{m\nu}\|_{L^3(\G)}^3dt
 \le\|u_0\|_H^2+C_0T+C_1TK_T+C_2\|f\|_{L^2(V^*)}^2.
\end{align*}

Let $A_0:V\to V^*$ and $A_1:V\to V^*$ be operators induced by $a_0$ and $a_1$. From \eqref{convbound} and the trace embedding,
\begin{align}\label{operatorbounds}
 \norm{A_0u_m}_{V^*}&\le C\norm{u_m}_V,\\
 \norm{A_1u_m}_{V^*}&\le C\norm{u_m}_V^2,\nonumber\\
 \norm{\gamma_\nu^*(|u_{m\nu}|u_{m\nu})}_{V^*}
 &\le C\norm{u_{m\nu}}_{L^3(\G)}^2,\nonumber\\
 \norm{\gamma_\nu^*g_m(t,u_{m\nu})}_{V^*}
 &\le C(1+\norm{u_m}_V).\nonumber
\end{align}
All terms on the right of the Galerkin equation are bounded in $L^1(0,T;V^*)$ by \eqref{apriori}. Becase in the Galerkin method we choose spacial kind of basis, from \eqref{proj} we have
\begin{equation}\label{deriv}
 u_m'=P_m\bigl(f-A_0u_m-A_1u_m
 -\gamma_\nu^*(\kappa|u_{m\nu}|u_{m\nu}+g_m(t,u_{m\nu}))\bigr),
 \qquad \norm{u_m'}_{L^1(0,T;V^*)}\le C_T.
\end{equation}

Let $F_m=f-A_0u_m-A_1u_m-
\gamma_\nu^*(\kappa|u_{m\nu}|u_{m\nu}+g_m(u_{m\nu}))$.
For each $v\in V$ with $\|v\|_V\le1$, H\"older's inequality yields
\begin{align*}
|a_1(u_m,u_m,v)|
&\le\|u_m\|_{L^4(\Om)}\|\nabla u_m\|_{L^2(\Om)}
\|v\|_{L^4(\Om)}\le C\|u_m\|_V^2,\\
\left|\int_{\G} |u_{m\nu}|u_{m\nu} v_\nu\,d\Ga\right|
&\le\|u_{m\nu}\|_{L^3(\G)}^2\|v_\nu\|_{L^3(\G)}
\le C\|r_m\|_{L^3(\G)}^2,\\
\left|\int_{\G} g_m(u_{m\nu})v_\nu\,d\Ga\right|
&\le\|g_m(u_{m\nu})\|_{L^2(\G)}\|v_\nu\|_{L^2(\G)}
\le C(1+\|u_m\|_V).
\end{align*}
Therefore
\begin{align*}
\int_0^T\|F_m(t)\|_{V^*}dt
&\le \sqrt T\|f\|_{L^2(V^*)}
 +C\sqrt T\|u_m\|_{L^2(V)}
 +C\|u_m\|_{L^2(V)}^2\\
&\quad+C T^{1/3}\|u_{m\nu}\|_{L^3(\Sigma_1)}^2
 +C(T+\sqrt T\|u_m\|_{L^2(V)})\le C_T.
\end{align*}
The uniform stability of $P_m$ in $V^*$ then implies
$\|u_m'\|_{L^1(V^*)}\le C_T$. Notice that the cubic boundary
term is controlled in $L^{3/2}$ in time, and the convection term
in $L^1$ in time.

From \eqref{apriori}--\eqref{deriv}, Simon's compactness theorem \cite{Simon} gives, after passage to a subsequence,
\begin{equation}\label{compactH}
 u_m\rightharpoonup u\ \text{in }L^2(0,T;V),\qquad
 u_m\to u\ \text{in }L^2(0,T;H).
\end{equation}

From Lemmma~\ref{lem:trace} ge get
\begin{equation}\label{traceconv}
 u_{m\nu}\to u_\nu\quad\text{strongly in }L^2(\Sigma_1),\quad
 u_{m\nu}(x,t)\to u_\nu(x,t)\quad\text{a.e. on }\Sigma_1.
\end{equation}
The $L^3$ bound also assures that $u_\nu\in L^3(\Sigma_1)$. Since $|u_{m\nu}|u_{m\nu}$ is bounded in $L^{3/2}(\Sigma_1)$ and converges pointwise almost everywhere, then
\begin{equation}\label{cubicweak}
 |u_{m\nu}|u_{m\nu}\rightharpoonup|u_\nu|u_\nu
 \quad\text{in }L^{3/2}(\Sigma_1).
\end{equation}
Similarly, \eqref{gm} and \eqref{traceconv} show that $g_m(\cdot,u_{m\nu})$ is bounded in $L^2(\Sigma_1)$; along a subsequence,
\begin{equation}\label{xiweak}
 g_m(\cdot,u_{m\nu})\rightharpoonup\xi\quad\text{in }L^2(\Sigma_1).
\end{equation}

\begin{lemma}\label{lem:graph}
Under $H_j^{\rm cub}$, convergences \eqref{traceconv} and \eqref{xiweak} imply
\begin{equation}\label{graph}
 \xi(x,t)\in\partial_r\psi(x,t,u_\nu(x,t))\quad\text{for almost every }(x,t)\in\Sigma_1.
\end{equation}
\end{lemma}
\begin{proof} Since the proof is through mollification, we give full proof, instead of just citing Aubin-Celina lemma.
 For a.e. $y=(x,t)\in \Sigma_1$, the multifunction
	$F_y(s)=\partial_r\psi(y,s)$ is upper semicontinuous with nonempty,
	compact and convex values. By the definition of the regularization,
	\begin{equation}\label{convhull}
		g_m(\cdot,u_{m\nu}(y))\in\overline{\operatorname{co}}
		\bigcup_{|s-u_{m\nu}(y)|\le 1/m}F_y(s).
	\end{equation}
	Indeed, the derivative of the mollified potential is a weighted average
	of the almost everywhere defined derivatives of $\psi(y,\cdot)$,
	which belong to its Clarke subdifferential.
	
	By \eqref{traceconv}, $u_{m\nu}(y)\to u_\nu(y)$ for a.e. $y\in S$.
	From the upper semicontinuity of $F_y$ and the convexity of its values,
	for every $\varepsilon>0$ and a.e. $y\in S$, there exists
	$m_0=m_0(y,\varepsilon)$ such that
	\[
	z_m(y)\in F_y(r(y))+\varepsilon[-1,1]
	\qquad\text{for all }m\ge m_0.
	\]
	Consequently,
	$\operatorname{dist}(g_m(\cdot,u_{m\nu}),F_y(u_\nu)\to0$ for a.e. $y\in S$.
	Moreover, by \eqref{xiweak} and the finiteness of the measure of $S$,
	$g_m(\cdot,u_{m\nu})\rightharpoonup\xi$ in $L^1(S)$. Hence, applying the
	Aubin--Cellina convergence theorem \cite{AubinCellina} to the
	measurable multifunction $y\mapsto F_y(u_\nu(y))$ with closed convex
	values, we obtain
	\[
	\xi(y)\in F_y(u_\nu(y))\quad\text{for a.e. }y\in S,
	\]
	which completes the proof.
\end{proof}

By \eqref{compactH} and weak continuity of $A_0$, the viscous term converges in the integrated weak formulation. For the convection term we can use compactness of $V\subset L^4(\Om;\R^2)$, together with \eqref{deriv}, to obtain
\begin{equation}\label{strongL4}
 u_m\to u\quad\text{strongly in }L^2(0,T;L^4(\Om;\R^2)).
\end{equation}
To see this directly, apply the same compactness criterion as above with $L^4(\Om)$ in place of $H$: the compact embedding $V\subset L^4(\Om)$ and the continuous injection $L^4(\Om)\hookrightarrow V^*$ are sufficient, and the time-translation estimate follows from \eqref{deriv} and the uniform $L^\infty(H)$ bound.

We pass to the limit with all terms. Take a vector $v\in V$. Put $v_m=P_mv$ and
$\varphi\in C_c^\infty(0,T)$. Since $v_m\to v$ in $V$,
$\|v_m\|_V\le C\|v\|_V$, and the convergence of the viscous
term follows from
\[
\int_0^T\varphi a_0(u_m,v_m)dt
=\int_0^T\varphi a_0(u_m,v)dt
 +\int_0^T\varphi a_0(u_m,v_m-v)dt.
\]
The first integral converges by weak convergence in $L^2(V)$.
The absolute value of the second is bounded by
\[
 C\|\varphi\|_\infty\|u_m\|_{L^2(V)}
 \sqrt T\|v_m-v\|_V\longrightarrow0.
\]
For the convection, write
\begin{align*}
a_1(u_m,u_m,v_m)-a_1(u,u,v)
&=a_1(u_m,u_m,v_m-v)\\
&\quad+a_1(u_m-u,u_m,v)
 +a_1(u,u_m-u,v).
\end{align*}
The first term, integrated in time, is bounded by
$C\|v_m-v\|_V\|u_m\|_{L^2(V)}^2$ and tends to zero.
For the second term, the inequality
\[
 |a_1(u_m-u,u_m,v)|
 \le\|u_m-u\|_{L^4(\Om)}\|\nabla u_m\|_{L^2(\Om)}
 \|v\|_{L^4(\Om)}
\]
gives
\[
\left|\int_0^T\varphi a_1(u_m-u,u_m,v)dt\right|
\le C\|\varphi\|_\infty\|v\|_V
\|u_m-u\|_{L^2(L^4)}\|u_m\|_{L^2(V)}\to0.
\]
For the third term, the function $w\mapsto
\int_0^T\varphi(t)a_1(u(t),w(t),v)dt$ is a bounded
linear functional on $L^2(0,T;V)$; indeed,
\[
\left|\int_0^T\varphi a_1(u,w,v)dt\right|
\le C\|\varphi\|_\infty\|v\|_V
\|u\|_{L^2(V)}\|w\|_{L^2(V)}.
\]
Thus weak convergence of $u_m-u$ in $L^2(V)$ applies.

The $L^3(\Sigma_1)$ estimate implies
\[
\||u_{m\nu}|u_{m\nu}\|_{L^{3/2}(\Sigma_1)}^{3/2}
=\int_{\Sigma_1}|u_{m\nu}|^3\le C_T.
\]
Strong convergence of traces in $L^2(\Sigma_1)$ gives, after
subsequence extraction, almost everywhere convergence. Therefore
$q_m\to |u_\nu|u_\nu$ almost everywhere. Boundedness in the reflexive
space $L^{3/2}$ and identification of all weak accumulation points
yield $u_{m\nu}|u_{m\nu}\rightharpoonup |u_\nu|u_\nu$ in $L^{3/2}$.
Since $v_{m\nu}\to v_\nu$ in $L^3(\G)$, we can split
\begin{align*}
\int_{\Sigma_1}\varphi(|u_{m\nu}|u_{m\nu}v_{m\nu}-|u_\nu|u_\nu v_\nu)
&=\int_{\Sigma_1}\varphi |u_{m\nu}|u_{m\nu}(v_{m\nu}-v_\nu)\\
&\quad+\int_{\Sigma_1}\varphi(|u_{m\nu}|u_{m\nu}-|u_\nu|u_\nu)v_\nu.
\end{align*}
The first integral tends to zero by H\"older's inequality and
uniform $L^{3/2}$ boundedness; the second by weak convergence.
The perturbation term is treated in exactly the same way, using
$g_m(u_{m\nu})\rightharpoonup\xi$ in $L^2(\Sigma_1)$ and
$v_{m\nu}\to v_\nu$ in $L^2(\G)$.

Finally, to recover the initial value without assuming $u'\in L^2(V^*)$,
fix $v\in V$ and integrate the Galerkin identity in time:
\begin{align*}
(u_m(t),P_mv)_H-(u_{0m},P_mv)_H
&=\int_0^t\langle f-A_0u_m-A_1u_m,v_m\rangle ds\\
&\quad-\int_0^t\int_{\G}
 (\kappa |u_{m\nu}|u_{m\nu}+g_m(u_{m\nu}))v_{m\nu}\,d\Ga ds.
\end{align*}
Every integrand is uniformly integrable in $L^1(0,T)$ as a scalar
function, though not necessarily uniformly equiintegrable for the
convection term. The distributional limit in time is sufficient:
for each test function $\varphi\in C_c^\infty(0,T)$ the preceding
weak convergences give the integrated equation for $u$.
The limiting right-hand side belongs to $L^1(0,T;V^*)$, so
$u$ has a $V^*$-absolutely continuous representative. The initial
value is obtained by testing the differentiated equation against
smooth time functions which need not vanish at $t=0$, using
$u_{0m}\to u_0$ in $H$. This yields $u(0)=u_0$ in $V^*$.
As $u\in L^\infty(0,T;H)$, the representative is weakly
continuous into $H$, and its initial value is $u_0$ in the weak
$H$ sense.

This completes the proof of Theorem~\ref{thm:main}.\qed

\section{Examples and related formulations}\label{sec:discussion}
\subsection{A nonsmooth, nonmonotone leak law}
The assumption $H_j^{\rm cub}$ is compatible with genuinely nonconvex and nonsmooth potentials. For instance, let $a,b\in L^\infty(\Sigma_1)$, with $a\ge0$, and define
\begin{equation}\label{examplej}
 \psi(x,t,r)=a(x,t)|r|-b(x,t)\cos r.
\end{equation}
Then
\[
 \partial_r\psi(x,t,r)=
 \begin{cases}
 \{a(x,t)\operatorname{sign}r+b(x,t)\sin r\},&r\ne0,\\
 [-a(x,t),a(x,t)],&r=0,
 \end{cases}
\]
and $|\partial_r\psi|\le\norm a_\infty+\norm b_\infty$. For suitable $b$, the potential is nonconvex. The total law
\[
 -\sigma_\nu\in\kappa|u_\nu|u_\nu+\partial_r\psi(x,t,u_\nu)
\]
therefore falls within Theorem~\ref{thm:main} for every $\kappa>1/2$.

A second example is obtained from a bounded measurable scalar function $g(x,t,\cdot)$ by setting $\psi(x,t,r)=\int_0^r g(x,t,s)\,ds$. Its Clarke subdifferential is contained in the interval between the essential lower and upper limits of $g$ near $r$. Thus discontinuous, nonmonotone functions are admitted as perturbations of the cubic law.

\subsection{The hemivariational inequality}
\begin{remark} The result can we written in the language of hemivariational inequalities, instead of inclusions. 	If $(u,\xi)$ is a weak solution, the defining property of the Clarke subdifferential gives, for every $v\in V$ and almost every $t$,
\[
 \int_{\G}\eta(t)v_\nu\,d\Ga\le
 \int_{\G}j^0(x,t,u_\nu(t);v_\nu)\,d\Ga.
\]
Consequently,
\begin{equation}\label{hvi}
 \du{u'(t)+A_0u(t)+A_1u(t)-f(t)}v
 +\int_{\G}j^0(x,t,u_\nu(t);v_\nu)\,d\Ga\ge0.
\end{equation}
\end{remark}

\begin{remark}
	
 The strict condition $\kappa>1/2$ is a sufficient condition ensuring a positive $L^3$ trace estimate; no optimality is claimed. The endpoint $\kappa=1/2$ does not provide the same trace coercivity through this argument.\end{remark}

\end{document}